\documentclass[11pt]{article}
\usepackage[centertags]{amsmath}
\usepackage{amsmath}
\usepackage{amssymb}
\usepackage{amsthm}
\usepackage{color}
\usepackage{latexsym}
\usepackage[T1]{fontenc}
\usepackage[latin1]{inputenc}
\usepackage{fancyhdr}
\usepackage[colorlinks=true,linkcolor=blue,citecolor=blue,urlcolor=blue]{hyperref} 
\numberwithin{equation}{section}

\numberwithin{equation}{section}
\newcommand{\be}{\begin{equation}}
\newcommand{\ee}{\end{equation}}

\newtheorem{definition}{Definition}[section]
\newtheorem{theorem}{Theorem}[section]
\newtheorem{proposition}{Proposition}[section]
\newtheorem{lemma}{Lemma}[section]

\newcommand{\ind}{\mathbf{1}}

\def \R{\mathbb{R}}

\def \N{\mathbb{N}}
\def \E{\mathbb{E}}

\def \bf{\textbf}

\def \bop {\noindent\textbf{Proof. }}
\usepackage{authblk}
\usepackage[english]{babel}
\usepackage{fancyhdr}

\begin{document}
\title{\huge{BSDEs with logarithmic growth driven by Brownian motion and Poisson random measure and connection to stochastic control problem}}
\author{{\Large{Khalid OUFDIL}}\\
{\Large e-mail: khalid.ofdil@gmail.com}\\Universit\'e Ibn Zohr, Equipe. Aide \'a la decision,
ENSA, B.P.  1136, \\Agadir, Maroc. }
\date{}
\maketitle

\noindent {\bf{Abstract.}}
In this paper, we study  one-dimensional backward stochastic differential equation  under logarithmic growth in the $z$-variable $(|z|\sqrt{|\ln|z|}|)$. We show the existence and the uniqueness of the solution when the noise is driven by a Brownian motion and an independent Poisson random measure. In addition, we highlight the connection of such BSDEs with stochastic optimal control problem, where we show the existence of an optimal strategy for the control problem. \\



\section{Introduction}
The notion of non-linear Backward Stochastic Differential Equations (BSDEs in short) was first introduced by Pardoux and Peng \cite{PP}. The solution of this equation, associated with a terminal value $\xi$ and a generator $f$, is a couple of stochastic processes $(Y_t,Z_t)_{t\leq T}$ such that
\begin{equation}
Y_t=\xi+\int_{t}^{T}f(s,Y_s,Z_s)ds-\int_{t}^{T}Z_sdB_s
\end{equation}
a.s. for all $t\leq T$, where $B$ is a Brownian motion and the processes  $(Y_t,Z_t)_{t\leq T}$ are adapted to a natural filtration of $B$. Since then, these equations have gradually became an important mathematical tool which is encountered in many fields such as financial mathematics, stochastic optimal control, partial differential equations, and other fields.
 
In this paper we are concerned with the problem of the existence and the uniqueness of a solution for one-dimensional BSDEs when the noise is driven by a Brownian motion and an independent Poisson random measure. Roughly speaking we look for a triple of adapted processes $(Y,Z,V)$ such that:
\begin{equation}
Y_t=\xi+\int_{t}^{T}f(s,Y_s,Z_s)ds-\int_{t}^{T}Z_sdB_s-\int_{t}^{T}\int_{\mathcal{U}}V_s(w)\tilde{\mu}(dw,ds)
\end{equation}
where $\tilde{\mu}$ is a compensated Poisson random measure.

Recently in \cite{P}, Kruse and Popier 
established the existence and the uniqueness
of a solution to multidimensional BSDEs in a general filtration which includes a Brownian motion and an independent Poisson random measure. The generator is under monotonic assumption w.r.t. the y-variable and the terminal condition and the driver are not necessarily square integrable. We however in this work focus on the existence and the uniqueness of the solution for BSDEs when the generator is under logarithmic growth assumption in the $z$-variable. It is worth noting that our technique is based on a localization method introduced in \cite{B1,B2} and more developed and extended in \cite{BEHP,BEH}.

It is well known that BSDEs are connected with stochastic optimal control. Therefore as an application of the obtained result, we show the existence of an optimal strategy for stochastic control of diffusion. Assume that we have a system whose dynamic is given by: 
 \begin{eqnarray}
 && x_t=x_0+\int_{0}^{t}\varphi(s,x_s,\bar{u}_s)ds+\int_{0}^{t}\sigma(s,x_s)dB^u_s+\int_{0}^{t}\int_{\mathcal{U}}\gamma(s,x_{s^-},w)\tilde{\mu}^u(dw,ds)\\\nonumber
 &&\qquad\qquad+\int_{0}^{t}\int_{\mathcal{U}}\gamma(s,x_{s^-},w)g(s,x_{s^-},\check{u}_s,w)\lambda(dw)ds.
 \end{eqnarray}
The agent intervenes in the system in the form of controlling with the help of the process $u=(\bar{u},\check{u})$. The idea is we characterize the value function as the unique solution of a specific BSDE under more relaxed assumptions. 

This paper is divided into five sections. In section 2, we begin by presenting our assumptions and the notations used trough out the paper. In section 3, we show some priori estimates for solutions of BSDEs. In section 4, we move on to study the existence and the uniqueness of the solution for the BSDEs using the localization method. Section 5 is devoted to the link between our BSDEs and the stochastic optimal control problem. We show   the value function as a solution of an appropriate BSDE.


\section{Notations and Assumptions}
\subsection{Notations}
Let $(\Omega,\mathcal{F},(\mathcal{F}_t)_{t\leq T},P)$ be a stochastic basis such that:
\begin{itemize}

\item $\mathcal{F}_0$ contains $\mathcal{N}$ the set of all $P$-null sets of $\mathcal{F}$,
\item $\mathcal{F}_{t^+}=\cap_{\epsilon>0}\mathcal{F}_{t+\epsilon}=\mathcal{F}_t,~~\forall t\leq T$.
\end{itemize}
   
We assume that $(\mathcal{F}_t)_{t\leq T}$ is supported by the two following independent processes:
\begin{enumerate}
\item[(i)] a standard $d$-dimensional Brownian motion $B=(B_t)_{0\leq t\leq T}$,
\item[(ii)] a Poisson random measure $\mu$ on $\mathcal{U}\times\R^+$ with intensity $\lambda(dw)dt$ where $\mathcal{U}\subset\R^m\setminus \{0\}$ $(m\geq 1)$. The corresponding compensated Poisson random measure $\tilde{\mu}(dw,dt)=\mu(dw,dt)-\lambda(dw)dt$ is a martingale w.r.t. $\mathcal{F}$. The measure $\lambda$ is $\sigma$-finite on $\mathcal{U}$ satisfying $$ \int_{\mathcal{U}}(1\wedge|w|^2)\lambda(dw)<+\infty.$$

\end{enumerate}
In this paper, let ${\mathcal{P}}$ denote the $\sigma$-algebra of $\mathcal{F}_t$-predictable sets on $\Omega\times[0,T]$. 

For a given adapted RCLL process $(X_t)_{t\leq T}$ and for any $t\leq T$ we set $X_{t^-}=\lim_{s\nearrow t}X_s$ with the convention that $X_{0^-}=X_0$ and $\Delta X_t=X_t-X_{t^-}$.
\newpage
Let us introduce the following spaces of processes and notations considered in this work, for
all $p \geq 1$
\begin{itemize}
\item $\mathcal{S}^p$ is the space of $\R $-valued $\mathcal{F}_t$-adapted and RCLL processes $\left(Y_t\right)_{t\in[0,T]}$ such that $$||Y||_{\mathcal{S}^p}=\E\left[\sup_{t\leq T}|Y_t|^p\right]^{\frac{1}{p}}<+\infty.$$
\item $\mathcal{M}^p$ denote the set of $\R^d$-valued $\mathcal{F}_t$-progressively measurable processes $\left(Z_t\right)_{t\in[0,T]}$ such that $$||Z||_{\mathcal{M}^p}=\E\left[\left(\int_{0}^{T}|Z_s|^2ds\right)^{\frac{p}{2}}\right]^{\frac{1}{p}}<+\infty.$$ 

\item $\mathcal{L}_{loc}$ is the space of all ${\mathcal{P}}\otimes\mathcal{B}(\mathcal{U})$-measurable mapping $V:\Omega\times[0,T]\times \mathcal{U}\rightarrow \R$ such that $\int_{0}^{T}\int_{\mathcal{U}}(|V_s(w)|^2\wedge|V_s(w)|)\lambda(dw)ds<+\infty$. Let $\mathcal{L}^p$ be the set of all $V\in\mathcal{L}_{loc}$ such that $$||V||_{\mathcal{L}^p}=\E\left[\left(\int_{0}^{T}\int_{\mathcal{U}}|V_s(w)|^2\lambda(dw)ds\right)^{\frac{p}{2}}\right]^{\frac{1}{p}}<+\infty.$$
\item $\mathcal{L}_{\lambda}^p=\mathcal{L}(\mathcal{U},\lambda;\R)$ is the set of measurable functions $V:\mathcal{U}\rightarrow \R$ such that
$$||V||^p_{\mathcal{L}_{\lambda}^p}=\int_{\mathcal{U}}|V(w)|^p\lambda(dw)<+\infty.$$  
\end{itemize}

\subsection{Assumptions}
Now let $\xi$ be and $\R$-valued and $\mathcal{F}_T$-measurable random variable and let $f:[0,T]\times\Omega\times\R\times\R^d\times\mathcal{L}_{\lambda}^2\rightarrow\R$  be a random function  which associates $(t,\omega,y,z,\nu)$ with $f(t,\omega,y,z,\nu)$. On the data $\xi$ and $f$ we make the following assumptions: 
\begin{enumerate}
\item[(\bf {H.1})] \ \ There exists a positive constant $A$ such that \ \ $\E \left[ |\xi|^{\ln(A T+2)+2}\right]<+\infty.$

\item[(\bf{H.2})]
\begin{enumerate}

\item[(i)]  $f$ is continuous in $(y,z,\nu)$ for almost all $(t,\omega)$;
    \item[(ii)] there exist two positives constants  $c_0$ and $c_1$ and a process $(\eta_t)_{t\leq T}$ such that for every $t, \omega, y, z, \nu$ :
        $$\mid f(t,\omega,y,z,\nu)\mid \leq |\eta_t|+c_0|z|\sqrt{\vert\ln(|z|)\vert}+c_1||\nu||_{\mathcal{L}^2_\lambda},$$

    \end{enumerate}
    where the process $(\eta_t)_{t\leq T}$ satisfies
    $$\E\left[\int_{0}^{T} |\eta_s|^{ \ln(Cs+2)+2}ds\right]<+\infty.$$
\item[(\bf {H.3})] \ There exist $ v^1:\Omega\times [0, T]\rightarrow\R^+$ \bigg(resp. $v^2:[0, T]\times\mathcal{U}\rightarrow\R^+$\bigg) satisfying $\E\left[\int_{0}^{T}({v^1}_s)^{q_1}ds\right]<+\infty$ \bigg(resp. $\E\left[\int_{0}^{T}||{v^2}_s||_{\mathcal{L}^2_\lambda}^{q_2}ds\right]<+\infty$\bigg) (for some $q_1$ and $q_2$)
and a constant $M \in\R^+$ such that: $\forall N>2,\; \hbox {and for every} \mid y\mid,\; \mid y'\mid,\; \mid z\mid,\; \mid z'\mid,\; ||\nu||_{\mathcal{L}^2_\lambda},\; ||\nu'||_{\mathcal{L}^2_\lambda}
 \leq \ln(N)$
 \begin{eqnarray*}
 &&\bar{y}
 \big(f(t,\omega,y,z,\nu)-f(t,\omega,y^{\prime},z',\nu')\big)
 \ind_{\{|{v^1}_t|+||{v^2}_t||_{\mathcal{L}^2_\lambda}\leq \ln(N)\}}\leq M|\bar{y}|^{2}\ln A_{N}\\
 &&\qquad\qquad\qquad\qquad\qquad\quad+ M\mid \bar{y}\mid\mid
 \bar{z}\mid\sqrt{\ln A_{N}}+M|\bar{y}|||\bar{\nu}||_{\mathcal{L}^2_{\lambda}}\sqrt{\ln A_N},
 \end{eqnarray*}
where $(\bar{y},\bar{z},\bar{\nu})=(y-y',z-z',\nu-\nu')$ and $(A_N)_{N>2}$ is a real valued sequence that satisfies  $1<A_N\leq \ln(N)^{r}$ for any $N>2$ and some $r>0$, moreover $\lim\limits_{N\rightarrow+\infty} A_N =+\infty.$
\end{enumerate}


\section{Estimations of solutions}
In this section we are going to provide estimates to the solution of the BSDE w.r.t. $\xi$ and $\eta$. However we first define the notion of solution of the  BSDE associated with the couple $(\xi,f)$.               
\begin{definition}
A triplet of processes $(Y,Z,V)=(Y_t,Z_t,V_t)_{t\leq T}$ is said to be a solution of the BSDE with jump associated with $(\xi,f)$ if the following holds:
\begin{equation}\label{equ}
\left\{\begin{array}{l}
Y\in\mathcal{S}^{\ln(A T+2)+2},~ Z\in\mathcal{M}^2 \mbox{ and }V\in\mathcal{L}^2;\\
Y_t =\xi + \int_t^T
f(s,Y_s,Z_s,V_s) ds- \int_t^T Z_s dB_s-\int_{t}^{T}\int_{\mathcal{U}}V_s(w)\tilde{\mu}(dw,ds),\qquad
t\in[0,T].
\end{array}\right.
\end{equation}
\end{definition}

The following lemma has already been mentioned and proved in \cite{BKK}, but for the sake of the reader we give it again.
\begin{lemma}\label{usyz}
Let $(Y_s,Z_s)\in\R\times \R^d$ and $s\in[0,T]$ be such that $Y_s$ is large enough. Then for every $C_2>0$ there exists a constant $C_3>0$ such that
\begin{equation}\label{majoration}C_2\mid
Y_s\mid|Z_s|\sqrt{\vert\ln(|Z_s|)\vert}\leq \frac{|Z_s|^2}{2}+ C_3\ln(\mid
Y_s\mid)\mid Y_s\mid^{2}.
\end{equation}
\end{lemma} 

Now we are ready to give estimates for $Y$, $Z$ and $V$. Actually we have the following proposition.
\begin{proposition}\label{estimate}
Let $(Y_t,Z_t,V_t)_{t\leq T}$ be a solution of the BSDE \eqref{equ}, and assume that the couple of data $(\xi,f)$ satisfies $(\bf {H.1})$ and $(\bf {H.2})$. Then there exists a constant $K>0$ such that:
\begin{equation}\label{estyv}
\E\left[\sup_{t\leq T}\mid Y_{t}\mid^{\theta(t)}+\int_0^{T} |Z_s|^2 ds+\int_{0}^{T}\int_{\mathcal{U}}|V_s(w)|^2\lambda(dw)ds\right]
\leq K\E\left[\mid
\xi\mid^{\theta(T)}+\int_{0}^{T}|\eta_s|^{\theta(s)}ds\right].
\end{equation}
where $\theta(t)=\ln(A t+2)+2$ for some constant $A$ large enough.
\end{proposition}
\bop: First let us show that for some positive constant $C$ we have
\begin{equation}\label{y}
\E\left[\sup_{t\leq T}\mid Y_{t}\mid^{\theta(t)}\right]
\leq C\E\left[\mid
\xi\mid^{\theta(T)}+\int_{0}^{T}|\eta_s|^{\theta(s)}ds\right].
\end{equation}
Let $\theta$ be the function from
$[0,T]$ into $\R^+$ defined by
$$\theta(t)=\ln(A t+2)+2,$$
and let $u$ be the function from $[0,T]\times\R$ into $\R^+$ defined by
$$u(t,x)=\mid x\mid^{\theta(t)}.$$
\newpage
Then
\begin{itemize}
\item $u_t=\frac{A}{A t+2}\ln(\mid x\mid)\mid x\mid^{\theta(t)}$;
\item $ u_x=\theta(t)\mid x\mid^{\theta(t)-2}x$;
\item $u_{xx}=\theta(t)(\theta(t)-1)\mid x\mid^{\theta(t)-2}$.
\end{itemize} 
Since there is a lack of integrability of the processes $Y$, $Z$ and $V$ we proceed by localization. Actually for $k\geq 0$, let
$\tau_{k}$ be the stopping time defined as follows:
$$\tau_k=\inf \{t\geq 0,\, \int_{0}^{t}\int_{\mathcal{U}}\theta(s)^2|Y_{s^-}|^{2\theta(s)-2}|V_s(w)|^2\lambda(dw)ds+\int_{0}^{t}\theta(s)^2\mid
Y_s\mid^{2\theta(s)-2}\mid Z_s\mid^2 ds\geq k\}\wedge T.$$
Using It\^o's formula on the process $Y_{t}$ and the function $y\rightarrowtail \mid y\mid^{\theta(t)}$ yields:
\begin{eqnarray}
&& \mid Y_{t\wedge \tau_k}\mid^{\theta(t\wedge\tau_k)} =\mid
Y_{\tau_k}\mid^{\theta(\tau_k)}- \int_{t\wedge \tau_k}^{\tau_k}\frac{A}{A s+2}\ln(\mid Y_s \mid)\mid  Y_s\mid^{\theta(s)}ds \\
\nonumber
&&\qquad\quad+\int_{t\wedge \tau_k}^{\tau_k}
\theta(s)\mid Y_s\mid^{\theta(s)-2}Y_sf(s,Y_s,Z_s,V_s)ds-\frac{1}{2}\int_{t\wedge
\tau_k}^{\tau_k} \theta(s)(\theta(s)-1)|Z_s|^2\mid
Y_s\mid^{\theta(s)-2}ds\\ \nonumber
&& \qquad~~ -  \int_{t\wedge
\tau_k}^{\tau_k}\theta(s)\mid
Y_s\mid^{\theta(s)-2}Y_sZ_s dB_s-\int_{t\wedge\tau_k}^{\tau_k}\int_{\mathcal{U}}\theta(s)Y_{s^-}|Y_{s^-}|^{\theta(s)-2}V_s(w)\tilde{\mu}(dw,ds)\\
\nonumber
&& \qquad~~ -\int_{t\wedge\tau_k}^{\tau_k}\int_{\mathcal{U}}\left(|Y_{s^-}+V_s(w)|^{\theta(s)}-|Y_{s^-}|^{\theta(s)}-\theta(s)Y_{s^-}|Y_{s^-}|^{\theta(s)-2}V_s(w)\right)\mu(dw,ds). \\ \nonumber
\end{eqnarray}
We use now the fact that $f$ satisfies (\bf {H.2}) and we obtain
\begin{eqnarray*}
&& \mid Y_{t\wedge \tau_k}\mid^{\theta(t\wedge\tau_k)} \leq \mid
Y_{\tau_k}\mid^{\theta(\tau_k)}- \int_{t\wedge \tau_k}^{\tau_k}\frac{A}{A s+2}\ln(\mid Y_s \mid)\mid  Y_s\mid^{\theta(s)}ds \\
\nonumber
&& \qquad \qquad+\int_{t\wedge \tau_k}^{\tau_k}
\theta(s)\mid Y_s\mid^{\theta(s)-1}(|\eta_s|+c_0|Z_s|\sqrt{\vert\ln(|Z_s|)\vert}+c_1||V_s||_{\mathcal{L}^2_\lambda})ds\\ \nonumber
&&\qquad\qquad-\frac{1}{2}\int_{t\wedge
\tau_k}^{\tau_k} \theta(s)(\theta(s)-1)|Z_s|^2\mid
Y_s\mid^{\theta(s)-2}ds-  \int_{t\wedge
\tau_k}^{\tau_k}\theta(s)\mid Y_s\mid^{\theta(s)-2}Y_s Z_s
dB_s\\ \nonumber
&& \qquad\qquad -\int_{t\wedge\tau_k}^{\tau_k}\int_{\mathcal{U}}\theta(s)Y_{s^-}|Y_{s^-}|^{\theta(s)-2}V_s(w)\tilde{\mu}(dw,ds)\\ \nonumber
&& \qquad\qquad -\int_{t\wedge\tau_k}^{\tau_k}\int_{\mathcal{U}}\left(|Y_{s^-}+V_s(w)|^{\theta(s)}-|Y_{s^-}|^{\theta(s)}-\theta(s)Y_{s^-}|Y_{s^-}|^{\theta(s)-2}V_s(w)\right)\mu(dw,ds). \\ 
\end{eqnarray*}
Next by Young's inequality it holds that 
$$
\theta(s)\mid Y_s\mid^{\theta(s)-1}|\eta_s| \leq \mid
Y_s\mid^{\theta(s)}+ \theta(s)^{\theta(s)-1}|\eta_s|^{
\theta(s)}.
$$
\newpage
\par\noindent
Now without loss of generality we assume that the $y$-variable is sufficiently large, then there exists a constant $C_1>0$ such that 
\begin{eqnarray*}
&& \mid Y_{t\wedge \tau_k}\mid^{\theta(t\wedge\tau_k)}\\
&&\qquad\qquad \leq \mid
Y_{\tau_k}\mid^{\theta(\tau_k)}+\int_{t\wedge \tau_k}^{\tau_k}\theta(s)^{\theta(s)-1}|\eta_s|^{\theta(s)}ds-\int_{t\wedge\tau_k}^{\tau_k}\int_{\mathcal{U}}\theta(s)Y_{s^-}|Y_{s^-}|^{\theta(s)-2}V_s(w)\tilde{\mu}(dw,ds)\\
\nonumber 
&& \qquad \qquad- \int_{t\wedge \tau_k}^{\tau_k}\theta(s)(\theta(s)-1)\mid Y_s\mid^{\theta(s)-2}\left[  \frac{C_1\ln(\mid Y_s \mid)\mid
Y_s\mid^{2}}{\theta(s)(\theta(s)-1)}+\frac{|Z_s|^2}{2}\right.\\ \nonumber
&&\qquad\qquad\left.-\frac{\theta(s)\mid
Y_s\mid c_0|Z_s|\sqrt{\vert\ln(|Z_s|)\vert}}{\theta(s)(\theta(s)-1)}\right]ds-  \int_{t\wedge
\tau_k}^{\tau_k}\theta(s)\mid Y_s\mid^{\theta(s)-2}Y_s Z_s
dB_s
\\\nonumber 
&& \qquad\qquad -\int_{t\wedge\tau_k}^{\tau_k}\int_{\mathcal{U}}\left(|Y_{s^-}+V_s(w)|^{\theta(s)}-|Y_{s^-}|^{\theta(s)}-\theta(s)Y_{s^-}|Y_{s^-}|^{\theta(s)-2}V_s(w)\right)\mu(dw,ds) \\ \nonumber 
&&\qquad\qquad+\int_{t\wedge \tau_k}^{\tau_k}
 c_1\theta(s)\mid Y_s\mid^{\theta(s)-1}||V_s||_{\mathcal{L}^2_\lambda}ds.
\end{eqnarray*}
Next there exist two positives constants $C_2$ and $C_3$ such that:

\begin{eqnarray}\label{fin-est}
&& \mid Y_{t\wedge \tau_k}\mid^{\theta(t\wedge\tau_k)} \leq \mid
Y_{\tau_k}\mid^{\theta(\tau_k)}+\int_{t\wedge \tau_k}^{\tau_k}\theta(s)^{\theta(s)-1}|\eta_s|^{ \theta(s)}ds- \int_{t\wedge \tau_k}^{\tau_k}\theta(s)(\theta(s)-1)\\
\nonumber && \qquad \qquad \mid Y_s\mid^{\theta(s)-2}\left[ C_3\ln(\mid Y_s
\mid)\mid Y_s\mid^{2}+\frac{|Z_s|^2}{2}-C_2\mid
Y_s\mid|Z_s|\sqrt{\vert\ln(|Z_s|)\vert}\right]ds
\\\nonumber 
&& \qquad\qquad -\int_{t\wedge\tau_k}^{\tau_k}\int_{\mathcal{U}}\left(|Y_{s^-}+V_s(w)|^{\theta(s)}-|Y_{s^-}|^{\theta(s)}-\theta(s)Y_{s^-}|Y_{s^-}|^{\theta(s)-2}V_s(w)\right)\mu(dw,ds) \\ \nonumber 
 && \qquad \qquad -  \int_{t\wedge
\tau_k}^{\tau_k}\theta(s)\mid Y_s\mid^{\theta(s)-2}Y_s Z_s
dB_s-\int_{t\wedge\tau_k}^{\tau_k}\int_{\mathcal{U}}\theta(s)Y_{s^-}|Y_{s^-}|^{\theta(s)-2}V_s(w)\tilde{\mu}(dw,ds)\\ \nonumber
&&\qquad\qquad+\int_{t\wedge \tau_k}^{\tau_k}
 c_1\theta(s)\mid Y_s\mid^{\theta(s)-1}||V_s||_{\mathcal{L}^2_\lambda}ds.
\end{eqnarray}
By using Lemma \ref{usyz} it follows that:
\begin{eqnarray}\label{fin-est1}
\nonumber
&& \mid Y_{t\wedge \tau_k}\mid^{\theta(t\wedge\tau_k)} \leq \mid
Y_{\tau_k}\mid^{\theta(\tau_k)}+\int_{t\wedge \tau_k}^{\tau_k}\theta(s)^{\theta(s)-1}|\eta_s|^{\theta(s)}ds+\int_{t\wedge \tau_k}^{\tau_k}
 c_1\theta(s)\mid Y_s\mid^{\theta(s)-1}||V_s||_{\mathcal{L}^2_\lambda}ds
\\
&& \qquad -\int_{t\wedge\tau_k}^{\tau_k}\int_{\mathcal{U}}\left(|Y_{s^-}+V_s(w)|^{\theta(s)}-|Y_{s^-}|^{\theta(s)}-\theta(s)Y_{s^-}|Y_{s^-}|^{\theta(s)-2}V_s(w)\right)\mu(dw,ds) \\ \nonumber 
 && \qquad -  \int_{t\wedge
\tau_k}^{\tau_k}\theta(s)\mid Y_s\mid^{\theta(s)-2}Y_s Z_s
dB_s-\int_{t\wedge\tau_k}^{\tau_k}\int_{\mathcal{U}}\theta(s)Y_{s^-}|Y_{s^-}|^{\theta(s)-2}V_s(w)\tilde{\mu}(dw,ds).
\end{eqnarray}
Next let us deal with the Poisson quantity
\begin{multline*}
-\int_{t\wedge\tau_k}^{\tau_k}\int_{\mathcal{U}}\left(|Y_{s^-}+V_s(w)|^{\theta(s)}-|Y_{s^-}|^{\theta(s)}-\theta(s)Y_{s^-}|Y_{s^-}|^{\theta(s)-2}V_s(w)\right)\mu(dw,ds)\\-\int_{t\wedge\tau_k}^{\tau_k}\int_{\mathcal{U}}\theta(s)Y_{s^-}|Y_{s^-}|^{\theta(s)-2}V_s(w)\tilde{\mu}(dw,ds).
\end{multline*}
Actually by following the same argument as in (\cite{P}, Proposition 2) it follows that 
\begin{eqnarray}\label{poisson}
\nonumber
&&-\int_{t\wedge\tau_k}^{\tau_k}\int_{\mathcal{U}}\left(|Y_{s^-}+V_s(w)|^{\theta(s)}-|Y_{s^-}|^{\theta(s)}-\theta(s)Y_{s^-}|Y_{s^-}|^{\theta(s)-2}V_s(w)\right)\mu(dw,ds)\\ 
&&\qquad\qquad-\int_{t\wedge\tau_k}^{\tau_k}\int_{\mathcal{U}}\theta(s)Y_{s^-}|Y_{s^-}|^{\theta(s)-2}V_s(w)\tilde{\mu}(dw,ds)\\ \nonumber
&&\qquad\qquad\leq -\int_{t\wedge\tau_k}^{\tau_k}\theta(s)(\theta(s)-1)3^{1-\theta(s)}|Y_{s^-}|^{\theta(s)-2}||V_s||_{\mathcal{L}^2_\lambda}^2ds\\\nonumber
&&\qquad\qquad-\int_{t\wedge\tau_k}^{\tau_k}\int_{\mathcal{U}}\left(|Y_{s^-}+V_s(w)|^{\theta(s)}-|Y_{s^-}|^{\theta(s)}\right)\tilde{\mu}(dw,ds).
\end{eqnarray}
We plug \eqref{poisson} in \eqref{fin-est1} and we obtain that: 
\begin{eqnarray}\label{fin-est2}
&& \mid Y_{t\wedge \tau_k}\mid^{\theta(t\wedge\tau_k)}+\int_{t\wedge\tau_k}^{\tau_k}\theta(s)(\theta(s)-1)3^{1-\theta(s)}|Y_{s^-}|^{\theta(s)-2}||V_s||_{\mathcal{L}^2_\lambda}^2ds
\\\nonumber 
&& \qquad\qquad\qquad \leq \mid
Y_{\tau_k}\mid^{\theta(\tau_k)}+\int_{t\wedge \tau_k}^{\tau_k}\theta(s)^{\theta(s)-1}|\eta_s|^{\theta(s)}ds-  \int_{t\wedge
\tau_k}^{\tau_k}\theta(s)\mid Y_s\mid^{\theta(s)-2}Y_s Z_s
dB_s\\ \nonumber
&& \qquad  \qquad\qquad-\int_{t\wedge\tau_k}^{\tau_k}\int_{\mathcal{U}}\left(|Y_{s^-}+V_s(w)|^{\theta(s)}-|Y_{s^-}|^{\theta(s)}\right)\tilde{\mu}(dw,ds)\\\nonumber
&&\qquad\qquad\qquad+\int_{t\wedge \tau_k}^{\tau_k}
c_1 \theta(s)\mid Y_{s^-}\mid^{\theta(s)-1}||V_s||_{\mathcal{L}^2_\lambda}ds.
\end{eqnarray}
Let $\alpha$ be a positive constant, once again by Young's inequality we have
$$
 \theta(s)\mid Y_{s^-}\mid^{\theta(s)-1}||V_s||_{\mathcal{L}^2_\lambda}\leq \frac{(c_1\theta(s))^2}{\alpha}|Y_s|^{\theta(s)}+\alpha |Y_{s^-}|^{\theta(s)-2}||V_s||^2_{\mathcal{L}^2_\lambda}.  $$
Therefore for $\alpha$ small enough
\begin{eqnarray}\label{fin-est3}
\nonumber
&& \mid Y_{t\wedge \tau_k}\mid^{\theta(t\wedge\tau_k)}\leq \mid
Y_{\tau_k}\mid^{\theta(\tau_k)}+\int_{0}^{\tau_k}\theta(s)^{\theta(s)-1}|\eta_s|^{\theta(s)}ds-  \int_{t\wedge
\tau_k}^{\tau_k}\theta(s)\mid Y_s\mid^{\theta(s)-2}Y_s Z_s
dB_s\\ \nonumber
&& \qquad \qquad+\int_{t\wedge\tau_k}^{\tau_k}\frac{(c_1\theta(s))^2}{\alpha}|Y_s|^{\theta(s)}ds-\int_{t\wedge\tau_k}^{\tau_k}\int_{\mathcal{U}}\left(|Y_{s^-}+V_s(w)|^{\theta(s)}-|Y_{s^-}|^{\theta(s)}\right)\tilde{\mu}(dw,ds).\\
&&
\end{eqnarray}
Then we take the expectation in \eqref{fin-est3} and since $\tau_k$ is of a stationary type we take the limit when $k\rightarrow+\infty$ and we obtain for some constant $C>0$ 
\begin{equation*}
\E\left[\mid Y_{t}\mid^{\theta(t)}\right]
\leq C\E\left[\mid
\xi\mid^{\theta(T)}+\int_{0}^{T}|\eta_s|^{\theta(s)}ds+\int_{0}^{T}|Y_s|^{\theta(s)}ds\right].
\end{equation*}
Finally we apply Gronwall's lemma and BDG inequality to get:
\begin{equation}\label{fin-est4}
\E\left[\sup_{t\leq T}\mid Y_{t}\mid^{\theta(t)}\right]
\leq C\E\left[\mid
\xi\mid^{\theta(T)}+\int_{0}^{T}|\eta_s|^{\theta(s)}ds\right].
\end{equation}
Now we will show that
\begin{multline*}
\E\left[\int_0^{T} |Z_s|^2 ds+\int_{0}^{T}\int_{\mathcal{U}}|V_s(w)|^2\lambda(dw)ds\right]\leq C \E\left[|\xi|^2 +
\sup\limits_{s\leq T}\mid Y_s\mid ^{\ln(2)+2} +
\int_0^T\mid \eta_s \mid^2ds\right].
\end{multline*}
Let $\delta_k$ be the sequence of stopping time defined as follow: $$\delta_k=\inf\{t\geq 0,~~\int_{0}^{t}|Z_s|^2ds+\int_{0}^{t}\int_{\mathcal{U}}|V_s(w)|^2\lambda(dw)ds\geq k\}\wedge T. $$\\
Applying once again It\^o's formula on the process $Y_t$ and the function $y\longmapsto y^2$ yields:
\begin{eqnarray*}
&& Y_0^2 + \int_0^{\delta_k} |Z_s|^2 ds+\int_{0}^{\delta_k}\int_{\mathcal{U}}|V_s(w)|^2\lambda(dw)ds = \xi^2 + 2\int_0^{\tau_k} Y_s
f(s,Y_s,Z_s,V_s) ds\\ \nonumber
&&\qquad\quad- 2 \int_0^{\delta_k}  Y_s Z_s dB_s-\int_{0}^{\delta_k}\int_{\mathcal{U}}\left(|Y_{s^-}+V_s(w)|^2-|Y_{s^-}|^2\right) \tilde{\mu}(dw,ds) \\ \nonumber
&& \qquad\quad\leq \xi^2 + 2\int_0^{\delta_k}\mid Y_s\mid(\mid
\eta_s\mid+c_0|Z_s|\sqrt{\vert\ln(|Z_s|)\vert}+c_1||V_s||_{\mathcal{L}^2_\lambda} )ds- 2 \int_0^{\delta_k} Y_s Z_s dB_s \\ \nonumber
&&\qquad\quad -\int_{0}^{\delta_k}\int_{\mathcal{U}}\left(|Y_{s^-}+V_s(w)|^2-|Y_{s^-}|^2\right) \tilde{\mu}(dw,ds).
\end{eqnarray*}

For any $\varepsilon>0$ we have
$
\sqrt{2\varepsilon \vert\ln(|z|)\vert}=\sqrt{\vert\ln(|z|^{2\varepsilon }\vert)}\leq
|z|^{\varepsilon }.$ Then by  using Young's inequality there exists a constant $C>0$ such that
\begin{eqnarray*}
&& \int_0^{\delta_k} |Z_s|^2 ds+\int_{0}^{\delta_k}\int_{\mathcal{U}}|V_s(w)|^2\lambda(dw)ds\leq |\xi|^2 + C\sup_{s\leq
T}|Y_s|^2+\int_0^T\mid \eta_s
\mid^2ds\\
&&\qquad\qquad~~+2\int_0^T \mid Y_s\mid (\frac{c_0}{\sqrt{2\varepsilon}}\mid Z_s\mid^{1+\varepsilon}+c_1||V_s||_{\mathcal{L}^2_\lambda} ) ds-2\int_0^{\delta_k} Y_s Z_s dB_s\\
&&\qquad\qquad~~-\int_{0}^{\delta_k}\int_{\mathcal{U}}\left(|Y_{s^-}+V_s(w)|^2-|Y_{s^-}|^2\right) \tilde{\mu}(dw,ds).
\end{eqnarray*}
Once again we use Young's inequality and we choose $0<\varepsilon <1$, then it holds
true that:
$$2c_0\frac{\mid Y_s\mid}{\sqrt{2\varepsilon}} \mid
Z_s\mid^{1+\varepsilon} \leq \frac{ 1-\varepsilon
}{2}(\frac{2c_0}{\sqrt{2\varepsilon}})^{\frac{2}{1-\varepsilon}}\mid Y_s\mid
^{\frac{2}{1-\varepsilon}}+\frac{1+\varepsilon}{2}\mid Z_s\mid^{2},$$
and for $\alpha>0$ we have: 
$$ 2c_1|Y_s|||V_s||_{\mathcal{L}^2_\lambda}\leq \frac{4c_1^2}{\alpha}|Y_s|^2+\alpha ||V_s||^2_{\mathcal{L}^2_\lambda}.$$
Then for $\varepsilon$ and $\alpha$ small enough, there exist positives constants $C_{\alpha}$, $C_\varepsilon$ and $c_{\varepsilon,\alpha}$ such that:
\begin{eqnarray*}
&& C_\varepsilon \int_0^{\delta_k} |Z_s|^2 ds+C_\alpha\int_{0}^{\delta_k}\int_{\mathcal{U}}|V_s(w)|^2\lambda(dw)ds \leq
|\xi|^2 +c_{\alpha,\varepsilon}\sup\limits_{s\leq T}\mid Y_s\mid^{\frac{2}{1-\varepsilon}}+\int_0^T\mid \eta_s
\mid^2ds\\
&&\qquad\qquad\qquad\qquad\qquad-2\int_0^{\delta_k} Y_s Z_s dB_s-\int_{0}^{\delta_k}\int_{\mathcal{U}}\left(|Y_{s^-}+V_s(w)|^2-|Y_{s^-}|^2\right) \tilde{\mu}(dw,ds).
\end{eqnarray*}
Now we take $\varepsilon \leq
\frac{\ln(2)}{2+\ln(2)}$, hence for $\mid y \mid$ large enough we have
\begin{eqnarray*}
\nonumber
&&  C_\epsilon \int_0^{\delta_k} |Z_s|^2 ds +C_\alpha\int_{0}^{\delta_k}\int_{\mathcal{U}}|V_s(w)|^2\lambda(dw)ds \leq |\xi|^2 +
c_{\alpha,\varepsilon}\sup\limits_{s\leq T}\mid Y_s\mid ^{\ln(2)+2} +
\int_0^T\mid \eta_s \mid^2ds\\ 
&&\qquad\qquad\qquad\qquad\qquad-2
\int_0^{\delta_k} Y_s Z_s dB_s-\int_{0}^{\delta_k}\int_{\mathcal{U}}\left(|Y_{s^-}+V_s(w)|^2-|Y_{s^-}|^2\right) \tilde{\mu}(dw,ds).\\ \nonumber
\end{eqnarray*} 
Therefor, there exists a positive constant $C$ such that 
\begin{eqnarray}\label{1}
\nonumber
&& \int_0^{\delta_k} |Z_s|^2 ds+\int_{0}^{\delta_k}\int_{\mathcal{U}}|V_s(w)|^2\lambda(dw)ds  \leq C\left(|\xi|^2 +
\sup\limits_{s\leq T}\mid Y_s\mid ^{\ln(2)+2} +
\int_0^T\mid \eta_s \mid^2ds\right.\\ \nonumber
&&\qquad\qquad\qquad\qquad\left.-2
\int_0^{\delta_k} Y_s Z_s dB_s-\int_{0}^{\delta_k}\int_{\mathcal{U}}\left(|Y_{s^-}+V_s(w)|^2-|Y_{s^-}|^2\right)\tilde{\mu}(dw,ds)\right).\\
&&
\end{eqnarray}
Next we take the expectation in \eqref{1} and we get 
\begin{multline*}
\E\left[\int_0^{\delta_k} |Z_s|^2 ds+\int_{0}^{\delta_k}\int_{\mathcal{U}}|V_s(w)|^2\lambda(dw)ds\right] \leq C \E\left[|\xi|^2 +
\sup\limits_{s\leq T}\mid Y_s\mid ^{\ln(2)+2} +
\int_0^T\mid \eta_s \mid^2ds\right].
\end{multline*}
Now since $\delta_k$ is a non decreasing sequence of stationary type which converge to $T$, we take the limit when $k\rightarrow+\infty$ and it follows that 
\begin{multline}\label{2}
\E\left[\int_0^{T} |Z_s|^2 ds+\int_{0}^{T}\int_{\mathcal{U}}|V_s(w)|^2\lambda(dw)ds\right] \leq C \E\left[|\xi|^2 +
\sup\limits_{s\leq T}\mid Y_s\mid ^{\ln(2)+2} +
\int_0^T\mid \eta_s \mid^2ds\right].
\end{multline}
Finally we combine \eqref{y} and \eqref{2} to complete the proof.

\section{Existence and uniqueness}
Before moving on to give the main theorem of this section we first introduce some  lemmas that will be useful in the proof of the existence and the uniqueness of the solution. 
\subsection{Useful Lemmas }
\begin{lemma}\label{estimf}
If (\textbf{H.2}) holds, then there exists a constant $c$ such that for $1<\bar{\alpha}<2$ we have
\begin{equation}
\E\left[\int_{0}^{T}|f(s,y,z,\nu)|^{\frac{2}{\bar{\alpha}}}ds\right]\leq c\E\left[1+\int_{0}^{T}|\eta_s|^2ds+\int_{0}^{T}|z|^2ds+\int_{0}^{T}||\nu||^2_{\mathcal{L}^2_\lambda}ds\right].
\end{equation}

\end{lemma}
\bop: From assumption (\textbf{H.2}) we can see that there exists a constant $\varepsilon>0$ such that 
\begin{equation*}
|f(t,y,z,\nu)|\leq |\eta_t|+\frac{c_0}{\sqrt{2\varepsilon}}|z|^{1+\varepsilon}+c_1||\nu||_{\mathcal{L}^2_\lambda}.
\end{equation*}
Then there exists a constant $c$ (that changes from line to line) such that 
\begin{eqnarray*}
&&|f(t,y,z,\nu)|^{\frac{2}{\bar{\alpha}}}\leq c\left(|\eta_t|^{\frac{2}{\bar{\alpha}}}+|z|^{\frac{2(1+\varepsilon)}{\bar{\alpha}}}+c_1||\nu||^{\frac{2}{\bar{\alpha}}}_{\mathcal{L}^2_\lambda}\right).
\end{eqnarray*}
Next we choose $0<\varepsilon<1$ and we put $\bar{\alpha}=1+\varepsilon$, then from Young's inequality we have 
\begin{equation*}
|f(t,y,z,\nu)|^{\frac{2}{\bar{\alpha}}}\leq c\left(1+|\eta_t|^2+|z|^2+||\nu||^2_{\mathcal{L}^2_\lambda}\right).
 \end{equation*}
 Finally we take expectation in both sides and it follows that 
 \begin{equation*}
\E\left[\int_{0}^{T}|f(s,y,z,\nu)|^{\frac{2}{\bar{\alpha}}}ds\right]\leq c\E\left[1+\int_{0}^{T}|\eta_s|^2ds+\int_{0}^{T}|z|^2ds+\int_{0}^{T}||\nu||^2_{\mathcal{L}^2_\lambda}ds\right].
 \end{equation*} \\
 
Now we will introduce a lemma whose proof is similar to the one in \cite{B}, and has an essential role in the proof of the existence of the solution.   
\begin{lemma}\label{exist}
 There exists a sequence of functions $(f_n)$ such that the following holds true:
 \begin{itemize}
 \item[{(i)}] For each $n$, $f_n$ is bounded and globally Lipschitz
 in $(y,z,\nu)$ $a.e.$ $t$ and $P$-$a.s.$ $\omega$.
 \item[(ii)] $ \forall n\geq 0$, $\mid f_n(t,\omega, y, z,\nu)\mid
 \leq |\eta_t|+c_0|z|\sqrt{\vert\ln(|z|)\vert}$, \quad $P$-$a.s.$, $a.e.$ $t\in
 [0,T]$.
 \item[(iii)] For every $N$,
 $\lim\limits_{n\rightarrow+\infty}\rho_N (f_n-f)=0$
 ; where $\rho_N(f)=E\left[\int_{0}^{T}\displaystyle\sup_{|y|,|z|,||\nu||_{\mathcal{L}_{\lambda}^2}\leq N}\mid f(s,y,z,\nu)\mid ds\right].$
 \end{itemize}  
\end{lemma}
\bop Let
$\psi_n$ from $\R^{3}$ to $ \R^+$  be a sequence of smooth functions such that:
\begin{equation*}
\psi_n(u)=\left\{\begin{array}{l}
1~~~\text{if}~~~\vert
u\vert \leq n;\\ 
0 ~~~\text{if} ~~~\vert u\vert > n.
\end{array}\right.
\end{equation*}
\noindent
We put, $\varepsilon_{q,n}(t,y,z,\nu) = \int
f(t,(y,z,\nu)-u)\alpha_q(u)du\psi_n(y,z,\nu)$;
where $\alpha_n: \R^3 \longrightarrow \R^+$ is a sequence of
smooth functions with compact support which approximate the Dirac
measure at 0 and which satisfy $\int \alpha_n (u)du = 1$. Next for $n \in \N^*$, let
$q(n)$ be an integer such that $q(n) \geq n+n^\alpha$. Then the sequence $f_n :=
\varepsilon_{q(n),n}$ satisfies all the assertions $(i)$-$(iii)$.\\

Next, Let $f_n$ be a sequence of functions associated to $f$ by Lemma \ref{exist}, and let us consider the following equation 
\begin{equation}\label{solf_n}
Y^n_t =\xi + \int_t^T
f_n(s,Y^n_s,Z^n_s,V^n_s) ds- \int_t^T Z^n_s dB_s-\int_{t}^{T}\int_{\mathcal{U}}V^n_s(w)\tilde{\mu}(dw,ds),\qquad
t\in[0,T].
\end{equation}
From Theorem 2.1 in \cite{BAR}, for each $n$ the Equation \eqref{solf_n} has a unique solution $(Y^n,Z^n,V^n)$. 
Then by using Proposition \ref{estimate} and Lemma \ref{estimf} we have the following lemma:

\begin{lemma}\label{Yn}
There exists a constant $K>0$ such that:
\begin{multline}
\E \left[\sup_{t\leq T}|Y^n_t|^{\ln(At+2)+2}+\int_0^{T} |Z^n_s|^2 ds+\int_{0}^{T}\int_{\mathcal{U}}|V^n_s(w)|^2\lambda(dw)ds\right.\\
\left.+\int_{0}^{T}|f_n(s,Y^n_s,Z^n_s,V^n_s)|^\frac{2}{\bar{\alpha}}ds\right]\leq K.
\end{multline}
\end{lemma}
 
 Now we first introduce the following lemma which has an essential role in the proof of the existence and uniqueness of the solution:
\begin{lemma}\label{Bet}
Let $\varepsilon>0$, $r>0$, $1<\bar{\alpha}<2$, $0<\kappa<2-\bar{\alpha}$ and  $2<\frac{2}{2-\bar{\alpha}}\leq \beta \leq \frac{(\ln(AT+2)+2)(2-\bar{\alpha}-\kappa)}{2}+1$. Assuming that (\textbf{H.1}), (\textbf{H.2}) and (\textbf{H.3}) are satisfied, let $(Y^n_t,Z^n_t,V^n_t)$ be the solution of \eqref{solf_n}. Then, there exists $K>0$ such that for appropriate $A$, $\bar{\alpha}$ and $\kappa$ and for $\delta<\frac{\kappa}{r(\beta M+2\beta M^2+\frac{3^{\beta-1}\beta M^2}{\beta-1})}$ we have:
\begin{equation}\label{estb}
 \lim\limits_{n,m\rightarrow+\infty}\E\left[\sup_{(T'-\delta)^+\leq t\leq T'}|Y^n_t-Y^m_{t}|^{\beta}\right] \leq \lim\limits_{n,m\rightarrow+\infty}K\E\left[e^{C_N\delta}|Y^n_{T'}-Y^m_{T'}|^{\beta}\right]+\varepsilon.
\end{equation}
\end{lemma}
\bop. We put $(\check{Y},\check{Z},\check{V})=(Y^n-Y^m,Z^n-Z^m,V^n-V^m)$, then by applying It\^o's formula on $e^{Ct}|\check{Y}_t|^\beta$ we obtain:
\begin{align*}
& e^{Ct}|\check{Y}_{t}|^{\beta}
+C\int_{t}^{T'}e^{Cs}|\check{Y}_{s}|^{\beta}ds\\
& =e^{CT'}|\check{Y}_{T'}|^{\beta}
+\beta\int_{t}^{T'}e^{Cs}|\check{Y}_{s}|^{\beta-2} \check{Y}_s
\big(f_{n}(s,Y_{s}^{n},Z_{s}^{n},V_s^n)-
f_{m}(s,Y_{s}^{m},Z_{s}^{m},V_s^m)\big)
ds
\\
& -\frac{\beta}{2}(\beta-1)\int_{t}^{T'}e^{Cs}|\check{Y}_{s}|^{\beta-2}\left|
\check{Z}_s\right|^{2}ds
-\beta\int_{t}^{T'}e^{Cs}|\check{Y}_{s}|^{\beta-2}
\check{Y}_s\check{Z}_s dB_{s}
\\
&-\beta\int_{t}^{T'}e^{Cs}\int_{\mathcal{U}}|\check{Y}_{s^-}|^{\beta-2}\check{Y}_{s^-} \check{V}_s(w)\tilde{\mu}(dw,ds)\\
&-\int_{t}^{T'}e^{Cs}\int_{\mathcal{U}}\left(|\check{Y}_{s^-}+\check{V}_s(w)|^{\beta}-|\check{Y}_{s^-}|^{\beta}-\beta|\check{Y}_{s^-}|^{\beta-2}\check{Y}_{s^-} \check{V}_s(w)\right)\mu(dw,ds).
\end{align*}
Put $\Phi(t)=|Y^n_t|+|Y^m_t|+|Z^n_t|+|Z^m_t|+||V^n_t||_{\mathcal{L}^2_{\lambda}}+||V^m_t||_{\mathcal{L}^2_{\lambda}}+|{v^1}_t|+||{v^2}_t||_{\mathcal{L}^2_{\lambda}}$, then 
\begin{align}\label{pp}
& e^{Ct}|\check{Y}_{t}|^{\beta}
+C\int_{t}^{T'}e^{Cs}|\check{Y}_{s}|^{\beta}ds+\frac{\beta}{2}(\beta-1)\int_{t}^{T'}e^{Cs}|\check{Y}_{s}|^{\beta-2}\left|
\check{Z}_s\right|^{2}ds\\\nonumber
& =e^{CT'}|\check{Y}_{T'}|^{\beta}
-\int_{t}^{T'}e^{Cs}\int_{\mathcal{U}}\left(|\check{Y}_{s^-}+\check{V}_s(w)|^{\beta}-|\check{Y}_{s^-}|^{\beta}-\beta|\check{Y}_{s^-}|^{\beta-2}\check{Y}_{s^-} \check{V}_s(w)\right)\mu(dw,ds)
\\\nonumber
& 
-\beta\int_{t}^{T'}e^{Cs}|\check{Y}_{s}|^{\beta-2}
\check{Y}_s\check{Z}_s dB_{s}-\beta\int_{t}^{T'}e^{Cs}\int_{\mathcal{U}}|\check{Y}_{s^-}|^{\beta-2}\check{Y}_{s^-}\check{V}_s(w)\tilde{\mu}(dw,ds)
\\\nonumber
&+J_{1}+J_{2}+J_{3}+J_{4},
\end{align}
where \begin{align*}
 & J_{1}:=\beta\int_{t}^{T'}e^{Cs}|\check{Y}_{s}|^{\beta-2} \check{Y}_s
\big(f_{n}(s,Y_{s}^{n},Z_{s}^{n},V_s^n)-
f_{m}(s,Y_{s}^{m},Z_{s}^{m},V_s^m)\big)
\ind_{\{\Phi(s)>\ln(N)\}} ds,
 \\ & J_{2}:=\beta\int_{t}^{T'}e^{Cs}|\check{Y}_{s}|^{\beta-2}
\check{Y}_s
\big(f_{n}(s,Y_{s}^{n},Z_{s}^{n},V_s^n)-
f(s,Y_{s}^{n},Z_{s}^{n},V_s^n)\big)
\ind_{\{\Phi(s)\leq \ln(N)\}} ds,
\\ & J_{3}:=\beta\int_{t}^{T'}e^{Cs}|\check{Y}_{s}|^{\beta-2}
\check{Y}_s
\big(f(s,Y_{s}^{n},Z_{s}^{n},V_s^n)-
f(s,Y_{s}^{m},Z_{s}^{m},V_s^m)\big)
\ind_{\{\Phi(s)\leq \ln(N)\}}ds,
\\ & J_{4}:=\beta\int_{t}^{T'}e^{Cs}|\check{Y}_{s}|^{\beta-2}
\check{Y}_s
\big(f(s,Y_{s}^{m},Z_{s}^{m},V_s^m)-
f_{m}(s,Y_{s}^{m},Z_{s}^{m},V_s^m)\big)
\ind_{\{\Phi(s)\leq \ln(N)\}}ds.
\end{align*}

\noindent
Now we will estimate $J_{1}$. let  $0<\kappa<2-\bar{\alpha}$, then we have 
\begin{align*}
J_{1} &\leq \beta e^{CT'} \dfrac{1}{\ln(N)^\kappa}\int_{t}^{T'}
|\check{Y}_{s}|^{\beta-1}{\Phi^\kappa(s)}
|f_{n}(s,Y_{s}^{n},Z_{s}^{n},V_s^n)-f_{m}(s,Y_{s}^{m},Z_{s}^{m},V_s^m)|ds
\\ &\leq
\beta e^{CT'} \dfrac{1}{\ln(N)^\kappa}\int_{t}^{T'}\left(|\check{Y}_{s}|^{\frac{2(\beta-1)}{2-\bar{\alpha}-\kappa}}+\Phi(s)^2+|f_{n}(s,Y_{s}^{n},Z_{s}^{n},V_s^n)-f_{m}(s,Y_{s}^{m},Z_{s}^{m},V_s^m)|^{\frac{2}{\bar{\alpha}}}\right)ds.
\end{align*}

\noindent
Next we estimate $J_{2}$ and $J_{4}$. It is easy to see that
\begin{align*}
J_{2}+ J_{4} & \leq 2\beta e^{CT'} (2\ln(N))^{\beta-1}
\bigg[\int_{t}^{T'}\sup_{|y|,|z|,||\nu||_{\mathcal{L}^2_{\lambda}}\leq
N}|f_{n}(s,y,z,\nu)-f(s,y,z,\nu)|ds
\\  &
+\int_{t}^{T'}\sup_{|y|,|z|,||\nu||_{\mathcal{L}^2_{\lambda}}\leq
N}|f_{m}(s,y,z,\nu)-f(s,y,z,\nu)|ds\bigg].
\end{align*}

\noindent
Finally we estimate $J_{3}$. Using assumption \bf{(H.3)}, we get
\begin{align*}
J_{3} \leq
\beta M\int_{t}^{T'}e^{Cs}
 \bigg[|\check{Y}_s|^{\beta}\ln A_{N}+|\check{Y}_s|^{\beta-1}|\check{Z}_s|\sqrt{\ln
A_{N}}+|\check{Y}_{s^-}|^{\beta-1}||\check{V}_s||_{\mathcal{L}^2_{\lambda}}\sqrt{\ln
A_{N}}  \bigg]ds.
\end{align*}
Now let us deal with the Poisson part in \eqref{pp}, as in Proposition 2 in \cite{P} we have 
\begin{eqnarray}
&&-\int_{t}^{T'}e^{Cs}\int_{\mathcal{U}}\left(|\check{Y}_{s^-}+\check{V}_s(w)|^{\beta}-|\check{Y}|_{s^-}^{\beta}-\beta|\check{Y}_{s^-}|^{\beta-2}\check{Y}_{s^-} \check{V}_s(w)\right)\mu(dw,ds)\\ \nonumber
&&\qquad\qquad\qquad-\beta\int_{t}^{T'}e^{Cs}\int_{\mathcal{U}}|\check{Y}_{s^-}|^{\beta-2} \check{Y}_{s^-} \check{V}_s(w)\tilde{\mu}(dw,ds)\\ \nonumber
&&\qquad\qquad\qquad\leq-\beta(\beta-1)3^{1-\beta}\int_{t}^{T'}e^{Cs}|\check{Y}_{s^-}|^{\beta-2}\int_{\mathcal{U}}|\check{V}_s(w)|^2\lambda(dw)ds\\ \nonumber
&&\qquad\qquad\qquad-\int_{t}^{T'}e^{Cs}\int_{\mathcal{U}}\left(|\check{Y}_{s^-}+\check{V}_s(w)|^{\beta}-|\check{Y}_{s^-}|^{\beta}\right)\tilde{\mu}(dw,ds).
\end{eqnarray}
Then \eqref{pp} becomes 

\begin{align*}
& e^{Ct}|\check{Y}_t|^{\beta}
+C\int_{t}^{T'}e^{Cs}|\check{Y}_s|^{\beta}ds+\frac{\beta}{2}\int_{t}^{T'}e^{Cs}|\check{Y}_s|^{\beta-2}\left|
\check{Z}_s\right|^{2}ds\\
& \leq e^{CT'}|\check{Y}_{T'}|^{\beta}
-\beta\int_{t}^{T'}e^{Cs}|\check{Y}_s|^{\beta-2}
\check{Y}_s\check{Z}_s dB_{s}-\int_{t}^{T'}e^{Cs}\int_{\mathcal{U}}\left(|\check{Y}_{s^-}+\check{V}_s(w)|^{\beta}-|\check{Y}_{s^-}|^{\beta}\right)\tilde{\mu}(dw,ds)
\\
&+
\beta e^{CT'} \dfrac{1}{\ln(N)^\kappa}\int_{t}^{T'}\left(|\check{Y}_{s}|^{\frac{2(\beta-1)}{2-\bar{\alpha}-\kappa}}+\Phi(s)^2+|f_{n}(s,Y_{s}^{n},Z_{s}^{n},V_s^n)-f_{m}(s,Y_{s}^{m},Z_{s}^{m},V_s^m)|^{\frac{2}{\bar{\alpha}}}\right)ds\\
&+\beta M\int_{t}^{T'}e^{Cs}
 \bigg[|\check{Y}_s|^{\beta}\ln A_{N}+|\check{Y}_s|^{\beta-1}|\check{Z}_s|\sqrt{\ln
A_{N}}+|\check{Y}_{s^-}|^{\beta-1}||\check{V}_s||_{\mathcal{L}^2_{\lambda}}\sqrt{\ln
A_{N}}  \bigg]\\
&+2\beta e^{CT'} (2\ln(N))^{\beta-1}
\bigg[\int_{t}^{T'}\sup_{|y|,|z|,||\nu||_{\mathcal{L}^2_{\lambda}}\leq
N}|f_{n}(s,y,z,\nu)-f(s,y,z,\nu)|ds
\\  &
+\int_{t}^{T'}\sup_{|y|,|z|,||\nu||_{\mathcal{L}^2_{\lambda}}\leq
N}|f_{m}(s,y,z,\nu)-f(s,y,z,\nu)|ds\bigg]\\
&-\beta(\beta-1)3^{1-\beta}\int_{t}^{T'}e^{Cs}|\check{Y}_{s^-}|^{\beta-2}\int_{\mathcal{U}}|\check{V}_s(w)|^2\lambda(dw)ds.
\end{align*}
Now let $\alpha'$ and $\gamma$ be two positives constants. By applying Young's inequality  we have 
\begin{eqnarray*}
&& \beta M|\check{Y}_{s}|^{\beta-1}|\check{Z}_s|\sqrt{\ln
A_{N}}\leq \frac{\beta^2 M^2\ln
A_{N}}{\alpha'}|\check{Y}_{s}|^{\beta}+\alpha'|\check{Y}_{s}|^{\beta-2}|\check{Z}_s|^2,
\end{eqnarray*}
and 
\begin{eqnarray*}
&&  \beta M|\check{Y}_{s^-}|^{\beta-1}||\check{V}_s||_{\mathcal{L}^2_{\lambda}}\sqrt{\ln
A_{N}}\leq \frac{\beta^2 M^2\ln
A_{N}}{\gamma} |\check{Y}_{s}|^{\beta}+\gamma|\check{Y}_{s^-}|^{\beta-2}||\check{V}_s||^2_{\mathcal{L}^2_{\lambda}}.
\end{eqnarray*}
Next we choose $\alpha'=\frac{\beta}{2}$, $\gamma=\beta(\beta-1)3^{1-\beta}$ and we put $C=C_N=\left(\beta M+2\beta M^2+\frac{3^{\beta-1}\beta M^2}{\beta-1}\right)\ln
A_{N}$, hence it follows
\begin{align*}
& e^{C_Nt}|\check{Y}_{t}|^{\beta} \leq 
-\beta\int_{t}^{T'}e^{C_Ns}|\check{Y}_{s}|^{\beta-2}
\check{Y}_s\check{Z}_s dB_{s}-\int_{t}^{T'}e^{Cs}\int_{\mathcal{U}}\left(|\check{Y}_{s^-}+\check{V}_s(w)|^{\beta}-|\check{Y}_{s^-}|^{\beta}\right)\tilde{\mu}(dw,ds)
\\
&\qquad+ \dfrac{\beta e^{C_NT'}}{\ln(N)^\kappa}\int_{t}^{T'}\left(|\check{Y}_{s}|^{\frac{2(\beta-1)}{2-\bar{\alpha}-\kappa}}+\Phi(s)^2+|f_{n}(s,Y_{s}^{n},Z_{s}^{n},V_s^n)-f_{m}(s,Y_{s}^{m},Z_{s}^{m},V_s^m)|^{\frac{2}{\bar{\alpha}}}\right)ds\\
&\qquad+e^{C_NT'}|\check{Y}_{T'}|^{\beta}+2\beta e^{C_NT'} (2\ln(N))^{\beta-1}
\bigg[\int_{t}^{T'}\sup_{|y|,|z|,||\nu||_{\mathcal{L}^2_{\lambda}}\leq
N}|f_{n}(s,y,z,\nu)-f(s,y,z,\nu)|ds
\\  &\qquad
+\int_{t}^{T'}\sup_{|y|,|z|,||\nu||_{\mathcal{L}^2_{\lambda}}\leq
N}|f_{m}(s,y,z,\nu)-f(s,y,z,\nu)|ds\bigg].
\end{align*}
We take the expectation and we use BDG inequality, there exists a constant $K>0$ such that 
\begin{align*}
& \E\left[\sup_{(T'-\delta)^+\leq t\leq T'}e^{C_Nt}|\check{Y}_{t}|^{\beta}\right] \leq K\E\left[e^{C_NT'}|\check{Y}_{T'}|^{\beta}\right]+ \dfrac{K\beta e^{C_NT'}}{\ln(N)^\kappa}\E\left[\int_{t}^{T'}|\check{Y}_s|^{\frac{2(\beta-1)}{2-\bar{\alpha}-\kappa}}ds\right]
\\
&\qquad\qquad+ \dfrac{K\beta e^{C_NT'}}{\ln(N)^\kappa}\E\left[\int_{t}^{T'}\left(\Phi(s)^2+|f_{n}(s,Y_{s}^{n},Z_{s}^{n},V_s^n)-f_{m}(s,Y_{s}^{m},Z_{s}^{m},V_s^m)|^{\frac{2}{\bar{\alpha}}}\right)ds\right]\\
&\qquad\qquad+2K\beta e^{C_NT'} (2\ln(N))^{\beta-1}
\E\bigg[\int_{t}^{T'}\sup_{|y|,|z|,||\nu||_{\mathcal{L}^2_{\lambda}}\leq
N}|f_{n}(s,y,z,\nu)-f(s,y,z,\nu)|ds
\\  &\qquad\qquad
+\int_{t}^{T'}\sup_{|y|,|z|,||\nu||_{\mathcal{L}^2_{\lambda}}\leq
N}|f_{m}(s,y,z,\nu)-f(s,y,z,\nu)|ds\bigg].
\end{align*}
Then it follows that
\begin{align*}
& \E\left[\sup_{(T'-\delta)^+\leq t\leq T'}|\check{Y}_{t}|^{\beta}\right] \leq K\E\left[e^{C_N\delta}|\check{Y}_{T'}|^{\beta}\right]
+\beta K \frac{A_N^{B\delta}}{A_N^{\frac{\kappa}{r}}}\E\left[\int_{t}^{T'}|\check{Y}_s|^{\frac{2(\beta-1)}{2-\bar{\alpha}-\kappa}}ds\right]\\
&\qquad\qquad+\beta K \frac{A_N^{B\delta}}{A_N^{\frac{\kappa}{r}}}\E\left[\int_{t}^{T'}\left(\Phi(s)^2+|f_{n}(s,Y_{s}^{n},Z_{s}^{n},V_s^n)-f_{m}(s,Y_{s}^{m},Z_{s}^{m},V_s^m)|^{\frac{2}{\bar{\alpha}}}\right)ds\right]\\
&\qquad\qquad+2K\beta e^{C_N\delta} (2\ln(N))^{\beta-1}
\E\bigg[\int_{t}^{T'}\sup_{|y|,|z|,||\nu||_{\mathcal{L}^2_{\lambda}}\leq
N}|f_{n}(s,y,z,\nu)-f(s,y,z,\nu)|ds
\\  &\qquad\qquad
+\int_{t}^{T'}\sup_{|y|,|z|,||\nu||_{\mathcal{L}^2_{\lambda}}\leq
N}|f_{m}(s,y,z,\nu)-f(s,y,z,\nu)|ds\bigg]
\end{align*}
where $B=\beta M+2\beta M^2+\frac{3^{\beta-1}\beta M^2}{\beta-1}$.
Then from  
Lemma \ref{exist} and Lemma \ref{Yn} we obtain for $\delta < \frac{\kappa}{Br}$ that 
\begin{equation}\label{estbb}
 \lim\limits_{n,m\rightarrow+\infty}\E\left[\sup_{(T'-\delta)^+\leq t\leq T'}|\check{Y}_{t}|^{\beta}\right] \leq \lim\limits_{n,m\rightarrow+\infty}K\E\left[e^{C_N\delta}|\check{Y}_{T'}|^{\beta}\right]+\varepsilon.
\end{equation}


\newpage
Now we are ready to give the main theorem of this section.
\begin{theorem}\label{ex}
 Assume that
\bf{(H.1)}, \bf{(H.2)} and \bf{(H.3)} are satisfied. Then the BSDE \eqref{equ} has a unique solution.
\end{theorem}
\bop. The proof will be divided into two parts, the first one is for the existence and the second one is for the uniqueness.\\

\noindent
\textbf{-Existence.}
Taking successively $T'=T$,
$T'=(T-\delta')^+$, $T'=(T-2\delta')^{+}...$ in \eqref{estb}
yields
\begin{equation}\label{esiYn}
\lim\limits_{n,m\rightarrow+\infty}\E\left[\sup_{0\leq t\leq T}|Y^n_{t}-Y^m_t|^{\beta}\right]=0.
\end{equation} 
Finally, Proposition \ref{estimate} allows us to show that
\begin{equation}\label{exiY}
\lim\limits_{n\rightarrow+\infty}\E\left[\sup_{0\leq t\leq T}|Y^n_t-Y_t|^{\beta}\right]=0.
\end{equation}
Next with $(\check{Y},\check{Z},\check{V})=(Y^n-Y^m,Z^n-Z^m,V^n-V^m)$, we apply once more It\^o's formula and we obtain:
\begin{eqnarray*}
&& \int_0^{T} |\check{Z}_s|^2 ds+\int_{0}^{T}\int_{\mathcal{U}}|\check{V}_s(w)|^2\lambda(dw)ds =  2\int_0^{T} \check{Y}_s
(f_n(s,Y^n_s,Z^n_s,V^n_s)-f_m(s,Y^m_s,Z^m_s,V^m_s)) ds\\ \nonumber
&&\qquad\quad- 2 \int_0^{T}  \check{Y}_s \check{Z}_s dB_s-\int_{0}^{T}\int_{\mathcal{U}}\left(|\check{Y}_{s^-}+\check{V}_s(w)|^2-|\check{Y}_{s^-}|^2\right) \tilde{\mu}(dw,ds)-|\check{Y}_0|^2. 
\end{eqnarray*}
We take the expectation, and by H\" older's inequality we have 
\begin{eqnarray*}
&& \E\left[\int_0^{T} |\check{Z}_s|^2 ds+\int_{0}^{T}\int_{\mathcal{U}}|\check{V}_s(w)|^2\lambda(dw)ds\right]\\
&&\qquad\qquad\leq 2\E\left[\int_0^{T} |\check{Y}_s|
|f_n(s,Y^n_s,Z^n_s,V^n_s)-f_m(s,Y^m_s,Z^m_s,V^m_s)| ds\right]\\ \nonumber
&&\qquad\qquad\leq C\E\left[ \sup_{0\leq t\leq T}|\check{Y}_s|^{2\over {2-\bar{\alpha}}}\right]^{\frac{2-\bar{\alpha}}{2}}\E\left[\int_0^{T} 
|f_n(s,Y^n_s,Z^n_s,V^n_s)-f_m(s,Y^m_s,Z^m_s,V^m_s)|^{2\over \bar{\alpha}}ds\right]^{\frac{\bar{\alpha}}{2}}.
\end{eqnarray*}
Then from \eqref{esiYn} and Lemma \ref{Yn} we have:
\begin{equation}
 \lim\limits_{n,m\rightarrow+\infty}\E\left[\int_0^{T} |Z^n_s-Z^m_s|^2 ds+\int_{0}^{T}\int_{\mathcal{U}}|V^n_s(w)-V^m_s(w)|^2\lambda(dw)ds\right]=0.
\end{equation}
Finally by using Proposition \ref{estimate} we get 
\begin{equation}
 \lim\limits_{n\rightarrow+\infty}\E\left[\int_0^{T} |Z^n_s-Z_s|^2 ds+\int_{0}^{T}\int_{\mathcal{U}}|V^n_s(w)-V_s(w)|^2\lambda(dw)ds\right]=0.
\end{equation}
Therefore, there exists a subsequence which we still denote $(Y^n,Z^n,V^n)$ such that:
\begin{equation}
\lim\limits_{n\rightarrow+\infty}\left(|Y^n_t-Y_t|+\int_0^{T} |Z^n_s-Z_s|^2 ds+\int_{0}^{T}\int_{\mathcal{U}}|V^n_s(w)-V_s(w)|^2\lambda(dw)ds\right)=0,~\text{a.e.}~ (t,\omega).
\end{equation}
On the other hand, as in the proof of Theorem 3.1 in \cite{B} we can show that: 
$$ \lim\limits_{n\rightarrow+\infty}\E\left[\int_{0}^{T}|f_n(s,Y^n_s,Z^n_s,V^n_s)-f(s,Y_s,Z_s,V_s)|ds\right]=0.$$
The existence is proved.\\

\noindent
\textbf{-Uniqueness.}\\

\noindent
Let $(Y,Z,V)$ and $(Y',Z',V')$ be two solutions of the BSDE \eqref{equ} associated with data $(\xi,f)$. Arguing as in  the proof of Lemma \ref{Bet} and by taking successively $T'=T$,
$T'=(T-\delta')^+$, $T'=(T-2\delta')^{+}...$, one can prove that $$Y_t=Y'_t~~~~~~~~~\forall t\in [0,T].$$ 
Next arguing once again as in the proof of the existence above, one can also prove that:

$$Z_t=Z_t' \text{ and }V_t=V'_t~~~~~~~~~\forall t\in [0,T].$$
The proof Theorem \ref{ex} is now complete.


\section{Application to stochastic optimal control problem}
Now we are going to highlight the link between stochastic optimal control, when the noise is of gaussian and Poisson types and BSDEs with Poisson jump.

Let $D$ be the space of controls, the space of the form $D_1\times D_2$, where $D_1$ and $D_2$ are two separable metric spaces whose Borel $\sigma$-algebra are $B(D_1)$ and $B(D_1)$ respectively.
Now let $\mathcal{D}_1$ (resp. $\mathcal{D}_2$) be the set of $\mathcal{P}$-measurable processes with values in $D_1$ (resp. $D_2$). $\mathcal{D}:=\mathcal{D}_1\times\mathcal{D}_2$ is called the set of admissible controls.\\

Let $x_0\in \R^d$ and $x_t$ be the solution of the following  standard functional differential equation.
\begin{equation}\label{x}
x_t=x_0+\int_{0}^{t}\sigma(s,x_s)dB_s+\int_{0}^{t}\int_{\mathcal{U}}\gamma(s,x_{s^-},w)\tilde{\mu}(dw,ds),
\end{equation}
where the mappings $\gamma: [0,T]\times \R^d\times \mathcal{U}\rightarrow \R^d$  and $\sigma:[0,T]\times \R^d\rightarrow\R^{d\times d}$ satisfy the following:
\begin{itemize}
\item[(i)] $\sigma$ is $\mathcal{P}$-measurable;
\item[(ii)] there exists a positive constant C such that $|\sigma(t,x)-\sigma(t,x')|\leq C|x-x'|$ and $|\sigma(t,x)|\leq C(1+|x|),$ for any $(x,x')$ in $\R^d\times \R^d$ and $t\leq T$;
\item[(iii)] for any $(t,x)\in [0,T]\times \R^d$, the matrix $\sigma(t,x)$ is invertible and $|\sigma^{-1}(t,x)|\leq C$ for some constant $C$;
\item[(iv)] for any $(t,x,x',w)\in [0,T]\times \R^d\times\R^d\times\mathcal{U}$, there exists a constant $C>0$ such that:
$$||\gamma(t,x,w)-\gamma(t,x',w)||_{\mathcal{L}^2_\lambda}\leq C|x-x'|~~~\text{and}~~~||\gamma(t,x,w)||_{\mathcal{L}^2_\lambda}\leq C(1+|x|)$$
\end{itemize}

According to Theorem 1.19 in \cite{OS} the process $(x_t)_{t\leq T}$ exists. Furthermore for $m\geq 2$ $(x_t)_{t\leq T}$ satisfies $\E\left[||x||_t^m\right]<+\infty$, for all $t$ where $||x||_t=\sup_{s\leq t}|x_s|$, (see \cite{BAR}, Proposition 1.1).

\newpage
Now we are given four functions $\varphi$, $c$, $g$ and $h$:
\begin{enumerate}
\item [(A.1)] $\varphi:[0,T]\times\R^d\times D_1\rightarrow\R^d$ such that: 
\begin{enumerate}
\item [(i)] for every $(t,x)\in[0,T]\times\R^d$ the function $\varphi(t,x,.): \bar{u}\rightarrow\varphi(t,x,\bar{u})$ is continuous,
\item [(ii)] there exists a real constant $K_1>0$ such that:
\begin{equation}
|\varphi(t,x,\bar{u})|\leq K_1(1+|x|) ~~~~\forall (t,x,\bar{u})\in[0,T]\times\R^d\times D_1;
\end{equation} 
\item[(iii)] there exists a positive constant $C$ such that 
\begin{equation}
|\sigma^{-1}(t,x)\varphi(t,x,\bar{u})|\leq C ~~~~\forall (t,x,\bar{u})\in[0,T]\times\R^d\times D_1;
\end{equation}
\end{enumerate}
\item [(A.2)] $c:[0,T]\times\R^d\times D_1\rightarrow\R$ such that: 
\begin{enumerate}
\item [(i)] for every $(t,x)\in[0,T]\times\R^d$ the function $c(t,x,.): \bar{u}\rightarrow c(t,x,\bar{u})$ is continuous, 
\item [(ii)]there exists a positive constant $K_2$ such that: 
\begin{equation}
|c(t,x,\bar{u})|\leq K_2(1+|x|) ~~~~\forall (t,x,\bar{u})\in[0,T]\times\R^d\times D_1;
\end{equation} 
\end{enumerate}
\item [(A.3)] $g:[0,T]\times\R^d\times D_2\times\mathcal{U}\rightarrow\R$ such that: 
\begin{enumerate}
\item [(i)] for every $(t,x,w)\in[0,T]\times\R^d\times\mathcal{U}$ the function $g(t,x,.,w): \check{u}\rightarrow g(t,x,\check{u},w)$ is continuous,
\item[(ii)] there exist two positives constants $\alpha_1$ and $\alpha_2$ such that for every $(t,x,\check{u},w)\in[0,T]\times\R^d\times D_2\times\mathcal{U}$
\begin{equation}
|g(t,x,\check{u},w)|\leq \alpha_1|w|\ind_{\{|w|\leq 1\}}+\alpha_2\ind_{\{|w|>1\}};
\end{equation}
\end{enumerate}
\item[(A.4)] $h:[0,T]\times\R^d\times D_2\times\mathcal{U}\rightarrow\R$ such that: 
\begin{enumerate}
\item[(i)] for every $(t,x,w)\in[0,T]\times\R^d\times\mathcal{U}$ the function $h(t,x,.,w): \check{u}\rightarrow h(t,x,\check{u},w)$ is continuous,
\item[(ii)]there exists a constant $C$ such that:  
\begin{equation}
||h(t,x,\check{u},w)||_{\mathcal{L}^2_\lambda}\leq C(1+|x|)~~~~\forall (t,x,\check{u},w)\in[0,T]\times\R^d\times D_2\times \mathcal{U};
\end{equation}
\end{enumerate}
\end{enumerate}

For any $u=(\bar{u},\check{u})\in \mathcal{D}$, let $L^u:=(L^u_t)_{t\leq T}$ be the positive local martingale solution of:  

$$dL^u_s=L^u_{s^-}\bigg\{\sigma^{-1}(s,x_s)\varphi(s,x_s,\bar{u}_{s})dB_s +\int_{\mathcal{U}}g(s,x_{s^-},\check{u}_{s},w)\tilde{\mu}(dw,ds) \bigg\} ~~\text{and}~~L^u_0=1.$$

Under the previous assumptions, for a given admissible control $u\in\mathcal{D}$ the process $L^u$ is an $(\mathcal{F}_t,P)$-martingale and $P^u$ defined by $dP^u=L^udP$ is a probability (\cite{BJ}, Corollary 5.1, pp. 244). Moreover under $P^u$, $\tilde{\mu}^u(dw,dt)=\tilde{\mu}(dw,dt)-g(t,x_{t^-},\check{u}_t,w)\lambda(dw)dt$ is an $\mathcal{F}_t$-martingale measure and $(B^u_t=B_t-\int_{0}^{t}\sigma^{-1}(s,x_s)\varphi(s,x_s,\bar{u}_s)ds)_{t\leq T}$  is an $\mathcal{F}_t$-Brownian motion.
 Then the Equation \eqref{x} becomes 
 \begin{eqnarray}
 && x_t=x_0+\int_{0}^{t}\varphi(s,x_s,\bar{u}_s)ds+\int_{0}^{t}\sigma(s,x_s)dB^u_s+\int_{0}^{t}\int_{\mathcal{U}}\gamma(s,x_{s^-},w)\tilde{\mu}^u(dw,ds)\\\nonumber
 &&\qquad\qquad+\int_{0}^{t}\int_{\mathcal{U}}\gamma(s,x_{s^-},w)g(s,x_{s^-},\check{u}_s,w)\lambda(dw)ds.
 \end{eqnarray}
 Now we shall study a stochastic control with one player. The controller, chooses an admissible control strategy $u^*$ to maximize the amount 
 \begin{equation}
 \int_{0}^{T} c(s,x_s,\bar{u}_s)ds+\int_{0}^{T}\int_{\mathcal{U}}h(s,x_{s^-},\check{u}_s,w)\lambda(dw)ds+\xi,
 \end{equation}   
  it is then in the interest of the controller to make the amount as big as possible at least on average, which lead us to a stochastic control with 
  \begin{equation}
  J(u)=\E^u\left[\int_{0}^{T} c(s,x_s,\bar{u}_s)ds+\int_{0}^{T}\int_{\mathcal{U}}h(s,x_{s^-},\check{u}_s,w)\lambda(dw)ds+\xi\right].
  \end{equation} 
But first let us introduce the two Hamiltonian functions $H_1$ and $H_2$ associated with this control problem defined on $[0,T]\times\R^d\times\R^d\times D_1$ and $[0,T]\times\R^d\times\mathcal{U}\times D_2$ respectively as follows 
\begin{equation}
H_1(t,x,z,\bar{u})=z\sigma^{-1}(t,x)\varphi(t,x,\bar{u})+c(t,x,\bar{u}).
\end{equation} 
and 
\begin{equation}
H_2(t,x,\nu,\check{u})=\int_{\mathcal{U}}\nu g(t,x,\check{u},w)\lambda(dw)+\int_{\mathcal{U}}h(t,x,\check{u},w)\lambda(dw).
\end{equation}
According to Benes selection theorem \cite{Be}, there exist two measurable functions $\bar{u}^*(s,x,z)$ and $\check{u}^*(s,x,\nu)$ with value in $D_1$ and $D_2$ such that 
$$H_1^*(t,x,z)=H_1(t,x,z,\bar{u}^*(s,x,z))=\sup_{\bar{u}\in D_1}H_1(t,x,z,\bar{u}),$$
and 
$$H_2^*(t,x,\nu)=H_2(t,x,z,\check{u}^*(s,x,\nu))=\sup_{\check{u}\in D_2}H_2(t,x,\nu,\check{u}).$$
Moreover the function $H^*(t,x,z,\nu)=H_1^*(t,x,z)+H_2^*(t,x,\nu)$ satisfies assumptions \textbf{(H.2)} and \textbf{(H.3)}. 
Indeed, there exists constants $C$, $c_0$ and $c_1$ such that:
$$|H^*_1(t,x,z)|\leq C(1+||x||_t)^2+c_0|z|\sqrt{|\ln|z||},$$
and 
$$|H^*_2(t,x,\nu)|\leq C(1+||x||_t)^2+c_1||\nu||_{\mathcal{L}^2_\lambda}.$$
Next to prove that the function $H^*$ satisfies \textbf{(H.3)} it is enough to take ${v^1}_t:=|\varphi(t,x,\bar{u})|^2$ and ${v^2}$ such that  $||{v^2}_t||_{\mathcal{L}^2_\lambda}:=||g(t,x,\check{u},w)||^2_{\mathcal{L}^2_\lambda}$.

Now we are ready to give the main result of this section. Actually we have the following theorem.

\begin{theorem} Assume that (A.1), (A.2), (A.3) and (A.4) are satisfied, and let $(Y^*,Z^*,V^*)$ be the solution of the BSDE associated with $(\xi,H^*)$. Then
the admissible control $u^*=(\bar{u}^*,\check{u}^*)$ is optimal strategy for the stochastic control; i.e. it satisfies
$$J(u^*)=Y^*_0\geq J(u)~~~~~~\forall u\in \mathcal{D}.$$
\end{theorem}

\bop Let $(Y^*,Z^*,V^*)$ be the solution of the BSDE associated with $(\xi,H^*)$. Then we have: 
 \begin{eqnarray*}
 && Y^*_0=\xi+\int_{0}^{T}H^*(s,x_s,Z^*_s,V^*_s)ds-\int_{0}^{T}Z^*_sdB_s-\int_{0}^{T}\int_{\mathcal{U}}V^*_s(w)\tilde{\mu}(dw,ds)\\
 &&\qquad=\xi+\int_{0}^{T}c(s,x_s,\bar{u}^*(s,x_s,Z^*_s))ds+\int_{0}^{T}\int_{\mathcal{U}}h(s,x_{s^-},\check{u}^*(s,x_s,V^*_s),w)\lambda(dw)ds\\
 &&\qquad-\int_{0}^{T}Z^*_sdB^{u^*}_s-\int_{0}^{T}\int_{\mathcal{U}}V^*_s(w)\tilde{\mu}^{u^*}(dw,ds).
 \end{eqnarray*}
Since $\left(\int_{0}^{t}Z^*_sdB^{u^*}_s\right)_{t\leq T} $ and $\left(\int_{0}^{t}\int_{\mathcal{U}}V^*_s(w)\tilde{\mu}^{u^*}(dw,ds)\right)_{t\leq T}$ are martingales, and since $Y^*_0$ is $\mathcal{F}_0$-measurable hence deterministic, then by taking the expectation we obtain 
 $$Y^*_0=\E^{u^*}\left[\xi+\int_{0}^{T}c(s,x_s,\bar{u}^*(s,x_s,Z^*_s))ds+\int_{0}^{T}\int_{\mathcal{U}}h(s,x_{s^-},\check{u}^*(s,x_s,V^*_s),w)\lambda(dw)ds\right].$$
It follows that
 $$Y^*_0=J(u^*).$$
 Let us now show that $Y^*_0\geq J(u)~~ \forall u\in\mathcal{D}.$ Actually
  \begin{eqnarray*}
  && Y^*_0=\xi+\int_{0}^{T}H^*(s,x_s,Z^*_s,V^*_s)ds-\int_{0}^{T}Z^*_sdB_s-\int_{0}^{T}\int_{\mathcal{U}}V^*_s(w)\tilde{\mu}(dw,ds)\\
  && \qquad\geq \xi+\int_{0}^{T}H(s,x_s,Z^*_s,V^*_s,\bar{u}_s,\check{u}_s)ds-\int_{0}^{T}Z^*_sdB_s-\int_{0}^{T}\int_{\mathcal{U}}V^*_s(w)\tilde{\mu}(dw,ds)\\
  &&\qquad= \xi+\int_{0}^{T}c(s,x_s,\bar{u}_s)ds+\int_{0}^{T}\int_{\mathcal{U}}h(s,x_{s^-},\check{u}_s,w)\lambda(dw)ds-\int_{0}^{T}Z^*_sdB^{u}_s\\
  &&\qquad-\int_{0}^{T}\int_{\mathcal{U}}V^*_s(w)\tilde{\mu}^{u}(dw,ds),
  \end{eqnarray*}
where $H=H_1+H_2$. Once more since $\left(\int_{0}^{t}Z^*_sdB^{u^*}_s\right)_{t\leq T} $ and $\left(\int_{0}^{t}\int_{\mathcal{U}}V^*_s(w)\tilde{\mu}^{u^*}(dw,ds)\right)_{t\leq T}$ are martingales we take the expectation and we get 
$$J(u^*)=Y^*_0\geq J(u).$$
This finishes the proof. \\

\noindent
\textbf{Funding:} This research was supported by National Center for Scientific and Technical Research (CNRST), Morocco.


\begin{thebibliography}{99}

\bibitem{B1}  K. Bahlali,   Backward stochastic differential
equations with locally Lipschitz coefficient, {C.R.A.S}, Paris,
serie I Math.  331, 481-486, 2001.

\bibitem{B2}  K. Bahlali,  Existence, uniqueness and stability for solutions of backward stochastic differential equations with locally Lipschitz coefficient, {Electron. Comm. Probab.}, 7, (2002), 169-179.

\bibitem{B} K. Bahlali and B. El Asri, Stochastic control and BSDEs with logarithmic growth, Bull. Sci. Math., 136(6) (2012), 617-637.

\bibitem{BEHP}  K. Bahlali, E. H. Essaky, M. Hassani and E. Pardoux, Existence, uniqueness and stability of backward stochastic differential equations with locally monotone coefficient, { C. R. Math. Acad. Sci.} Paris 335, (2002), no. 9, 757-762.

\bibitem{BEH} K. Bahlali, E. H. Essaky and M. Hassani,  Multidimensional BSDEs with super-linear growth coefficient: application to degenerate systems of semilinear PDEs, {C. R. Math. Acad. Sci.} Paris 348, 2010, no. 11-12, 677-682.


\bibitem{BKK} K. Bahlali, O. Kebiri, N. Khelfallah and H. Moussaoui, One dimensional BSDEs with logarithmic growth application to PDEs., stochastics, 2017. 

\bibitem{BAR} G. Barles, R. Buckdahn and E. Pardoux, BSDEs and integral-partial differential equations, Stochastics 60, 57-83, 1997.

\bibitem{Be} V. E. Benes, Existence of optimal stochastic control laws, SIAM JCO, 9, (3) 446-472, 1971.

\bibitem{BJ} A. Bensoussan and J. L. Lions, {Contr\^ole impulsionnel et in\'equations quasi-variationnelles}, Dunod, Paris, 1982.

\bibitem{FK} { T. Fujiwara and H. Kunitha}, { Stochastic differential equations of Jump type and  L\'{e}vy processes in diffeomorphism group}, J. Math. Kyoto  Univ., 25(1), 71-106, 1989.

\bibitem{HO} S. Hamad\`ene and Y. Ouknine, Reflected backward stochastic differential equations with jumps and random obstacles, Electronic Journal of Probability Vol. 8, 2003.

\bibitem{HW} S. Hamad\`ene and H. Wang, BSDEs with two RCLL Reflecting Obstacles driven by a Brownian Motion and Poisson Measure and related Mixed Zero-Sum Games, Stochastic Processes and their Applications, vol. 119, pp. 2881-2912, 2009. 

\bibitem{P} { T. Kruse and A. Popier}, { BSDEs with jumps in a general filtration}, Stoch. Int. J. Probab. Stoch. Proces, 88, Vol(2016): p491-p539.

\bibitem{OS} { B. \O ksendal and A. Sulem}, {Applied Stochastic Control of Jump Diffusion}, Springer Universitext, 2005.

\bibitem{LP} { J. P. Lepeltier and B.Marchal}, { Existence de politique optimale dans le contr\^ole int\'egro-diff\'erentiel},
Annales de l'I.H.P., vol. XIII, n.1, 1977, p.45-97.

\bibitem{PP} { E. Pardoux and S. Peng}, { Adapted solution of a
backward stochastic differential equation}, System Control Lett., \bf{14}, pages 55-61, 1990.

\bibitem{Y} { S. Yao}, { $L^p $  Solutions of Backward Stochastic Differential Equations with jumps}, (July 7 2016).


\end{thebibliography}
\end{document}